\documentclass[english,12pt]{article}

\usepackage{amsmath}
\usepackage{amssymb}
\usepackage{graphicx}
\usepackage{epic,eepic}

\usepackage{color}

\usepackage{amstext,amsthm,amssymb,amsmath}
\usepackage{graphicx}
\usepackage{subfigure}
\usepackage{wrapfig}
\usepackage{fullpage}
\usepackage{color}
\usepackage{multirow}
\usepackage{tabulary}
\usepackage{booktabs}
\usepackage{enumerate}
\usepackage{eepic,epic}
\usepackage{epsfig,subfigure,epstopdf}
\usepackage{algorithm,algorithmicx,algpseudocode}

\newtheorem{theorem}{Theorem}[section]

\newtheorem{remark}[theorem]{Remark}

\hfuzz=\maxdimen
\title{Cluster-based Generalized Multiscale Finite Element Method for elliptic PDEs with random coefficients}

\author{
Eric T. Chung \thanks{Department of Mathematics,
The Chinese University of Hong Kong (CUHK), Hong Kong SAR. Email: {\tt tschung@math.cuhk.edu.hk}.}
\and
Yalchin Efendiev \thanks{Department of Mathematics \& Institute for Scientific Computation (ISC),
Texas A\&M University,
College Station, TX, 77843, USA. Email: {\tt efendiev@math.tamu.edu}.}
\and
Wing Tat Leung \thanks{Department of Mathematics, Texas A\&M University, College Station, TX 77843, USA.
Email: {\tt sidnet123@yahoo.com.hk}.}
\and
Zhiwen Zhang \thanks{Department of Mathematics, The University of Hong Kong, Hong Kong SAR.
Email: {\tt zhangzw@hku.hk}. Corresponding author. }
}

\date{}

\usepackage{babel}
\begin{document}

\maketitle

\begin{abstract}
We propose a generalized multiscale finite element method (GMsFEM) based on clustering algorithm to study the elliptic PDEs with random coefficients in the multi-query setting.  Our method consists of offline and online stages. In the offline stage, we construct a small number of reduced basis functions within each coarse grid block, which can then be used to approximate the multiscale finite element basis functions. In addition, we coarsen the corresponding random space through a clustering algorithm. In the online stage, we can obtain the multiscale finite element basis very efficiently on a coarse grid by using the pre-computed multiscale basis. The new GMsFEM can be applied to multiscale SPDE starting with a relatively coarse grid, without requiring the coarsest grid to resolve the smallest-scale of the solution.  The new method offers considerable savings in solving multiscale SPDEs. Numerical results are presented to demonstrate the accuracy and efficiency of the proposed method for several multiscale stochastic problems without scale separation. 
\end{abstract}

\section{Introduction}

Many multiscale problems have stochastic nature due to missing information at small scales. For example, in porous media applications, the media properties are unknown in many locations and interpolated from some indirect and direct measurements. As a result, engineers can typically produce many realizations of the permeability field. In some cases, the stochastic description of the permeability is parameterized; however, in many cases, one deals with very large number of permeability realizations. The objective is to perform
many simulations and understand the solution as a stochastic quantity. In this paper, our objective is to propose a fast method based  on coarsening of both spatial and uncertainty space (in terms of realizations) for computing the solution space.

The stochastic description of the permeability field contains multiple scales, which do not have apparent scale
separation. Moreover, the uncertainties at different scales can be tightly coupled and one needs to upscale their interaction together. For this reason, stochastic upscaled models are often used. Stochastic upscaled models use coarse grids in the physical space and propose an ensemble average of the fine-grid permeability. As a result, one deals with stochastic flow equations on coarse grids. Some more commonly used approaches include upscaling a few realizations and making predictions based on these simulations.

Many  approaches have been developed for stochastic problems.  Most famous approaches include Monte Carlo methods and their variations \cite{glasserman:03,Caflisch:98,Niederreiter:01,Giles:08}. In some other approaches the uncertainties are represented using polynomials of random fields or collocation points in uncertainty space. Typical methods are the stochastic finite element method [2], Wiener chaos expansion or generalized polynomial chaos (gPC) method and stochastic collocation method, see \cite{Ghanem:91,Knio:01,Xiu:03,babuska:04,matthies:05,WuanHou:06,Wan:06,Xiu:09,Najm:09,sapsis:09} and references therein. These approaches do not take into account the interaction between the spatial scales and uncertainties. There are some recent attempts
that take into account of the space and uncertainty interaction
\cite{nouy2007generalized,nouy2010proper}.
In \cite{Zabaras:13}, Zabaras et al proposed a probabilistic graphical model approach to efficiently perform uncertainty quantification in the presence of both stochastic input and multiple physical scales. Hou et al explored the low-dimensional structures that hidden behind the high-dimensional SPDE problems and
have made some progress in solving high-dimensional and/or multiscale SPDEs \cite{ChengHouZhang1:13,ChengHouZhang2:13,ChengHouYanZhang:13,ZhangCiHouMMS:15}. We should point out that it is still extremely challenging to solve stochastic multiscale problems, due to their complicated natures,
namely, high-dimensionality of the solution spaces.


In our approach, we consider the joint coarsening of both spatial and uncertainty space motivated by many physical applications,
where uncertainties and spatial scales are tightly coupled. For representing the spatial scales, we use the generalized multiscale
finite element method (GMsFEM)
\cite{Efendiev:13,chung2016adaptive,egh12,ge09_2,hw97,ehw99}.
The main idea of this approach is to use the snapshots and local spectral decomposition for computing few basis functions. In this paper, these ideas are used for a few selected realizations to compute basis
functions that can be used to represent multiple realizations. However, it is obvious that the use of many realizations will introduce many spatial patterns and can require too many multiscale basis functions. This can lead to very expensive computations.  For this reason, we use coarsening algorithm in the uncertainty space that can allow using a few multiscale basis functions and possibly basis functions that can be written as a product of spatial and stochastic functions. In general, if the uncertainty space is partitioned uniformly, we
expect a larger dimensional coarse spaces.

To coarsen the uncertainty space, we first define a distance function in each coarse block. The appropriate distance functions
in clustering must include the distance between the solutions as this is a correct measure. However, one usually has many
local solutions. In this case, we use local solutions with randomized boundary conditions motivated by \cite{randomized2014}. The main idea of
our motivation is that the local problems with randomized boundary conditions can represent the main modes of local problems
that are computed using all boundary conditions. Because there are many realizations, we compute these local
solutions for selected realizations and then compute their Karhunen-Lo\`{e}ve expansion (KLE) \cite{Karhunen:47,Loeve:78} to approximate the solutions at all other realizations.  Using these approximate solutions, we compute the distance between two realizations. Using this distance and k-mean clustering algorithms \cite{Hartigan:1979}, we cluster the uncertainty space and define the corresponding coarse blocks.

To construct multiscale basis functions in each cluster and a coarse block, we select a realization or an average realization. The latter is used within GMsFEM to construct multiscale basis functions. The basis function is constructed by computing snapshot vectors and performing local spectral decomposition. In this construction, we assume that the permeability fields within each cluster is almost uniform in the appropriate metric as defined for clustering. One can also try to construct product basis functions. To do so,
we will need to construct stochastic basis functions, which represent the variations of the solution in the uncertainty space. Taking the product of stochastic basis functions and spatial multiscale basis functions for some average realizations, we can compute
multiscale basis functions, which represent both spatial information and uncertainties.

For global coupling of multiscale basis functions, we use continuous Galerkin formulation, though other discretization approaches can be used (e.g., Petrov-Galerkin, Discontinuous Galerkin methods). We consider two global couplings. The first one is a realization based approach, where we solve the global problem for each realization. The computations for each realization can be performed in parallel. Another approach is to use Galerkin in space and uncertainty space. In this approach, the variational formulation takes into account the average across the realizations. The second approach is cheaper for ensemble level calculations because it involves fewer global computations.

The offline space constructed in this paper can be enhanced using online basis functions. The online basis functions are adaptively constructed using the residual information to improve the convergence. In \cite{chung2015residual}, this algorithm is discussed for deterministic problems. The extension of this approach to stochastic problems will be discussed in this paper. The online basis functions are constructed for spatial and stochastic spaces by considering the joint residual. The construction can be done adaptively for selected coarse blocks.

The rest of the paper is organized as follows. The formulations of the stochastic multiscale model problem is presented in
section 2. In section 3, we introduce the coarse grid discretization of the spacial domain and the random space.
We also propose the basic clustering algorithm for coarsening the tensor space.
In section 4, we give a general introduction of the generalized multiscale finite element method (GMsFEM) based on clustering. Both the offline and online computation will be discussed. In addition, issues regarding the practical implementation and computational complexity analysis will be covered. In Section 5, we present numerical results to demonstrate the efficiency of our method. Concluding remarks are made in Section 6.


\section{Preliminaries}

In this section, we present the problem setting and the fine grid
discretization. Let $D$ be a bounded domain in
$\mathbb{R}^{d}\, (d=2,3)$
and $\Omega$ is a parameter space in $\mathbb{R}^{m}$, where the dimension $m$ can be large in practice.
We consider
the following parameter-dependent second order linear elliptic equation,
\begin{equation}
\label{eq:PDE}
\begin{cases}
-\nabla\cdot(\kappa(x,\omega)\nabla u(x,\omega))=f(x,\omega), & \text{ }(x,\omega)\in D\times\Omega, \\
u(x,\omega)=0, &\text{ }(x,\omega)\in \partial D\times\Omega,
\end{cases}
\end{equation}

where $\kappa(x,\omega)$ is highly heterogeneous with respect to the physical space $D$ and $f(x,\omega)\in L^2(D)$,
$\forall \omega \in \Omega$. Even though we write the problem as parametric, our method will only use
realizations of $\kappa(x,\omega)$. This is very important since engineers can only obtain realizations of the permeability field in practice.
Thus, we assume that $\Omega$ is a discrete set of realizations $\omega_i$'s.
The weak formulation
of the problem (\ref{eq:PDE}) is to find $u(x,\omega)\in L^{2}(\Omega;H^{1}(D))$ such that
\begin{equation}
\label{eq:weak}
\int_{\Omega} a(\omega;u,v) \; d\omega=\int_{\Omega}l(\omega;v) \; d\omega, \quad\quad \forall v\in L^{2}(\Omega;H^{1}(D)),
\end{equation}
where the bilinear form $a$ and the linear functional $l$ are defined as,
\begin{align*}
a(\omega;u,v) & =\int_{D}\kappa(x,\omega)\nabla u(x,\omega)\cdot\nabla v(x,\omega)dx,\\
l(\omega;v) & =\int_{D}f(x,\omega)v(x,\omega)dx.
\end{align*}
We remark that the measure $d\omega$ in (\ref{eq:weak}) is assumed to be the uniform measure as in many cases we only have realizations of the permeability field $\kappa(x,\omega)$ (e.g., the SPE10 model in reservoir simulation). In the space $L^{2}(\Omega;H^{1}(D))$, we define the norm of the tensor space $H^{1}_{0}(D)\otimes L^{2}(\Omega)$ as,
\begin{equation}
\| u\|_{H^{1}_{0}(D)\otimes L^{2}(\Omega)}^2 = \int_{\Omega} a(\omega;u,u) \; d\omega.
\end{equation}

\subsection{Fine grid discretization }

We introduce the fine grid discretization for solving the multiscale
problem (\ref{eq:weak}). Let $\mathcal{T}^{h}$ be a fine-grid partition of the domain
$D$ with mesh size $h$ and $\Omega_{d}$ be the fine sampling subset
of $\Omega$ which can approximately represent the whole space $\Omega$,
namely $\Omega_{d}=\{\omega_{1},\omega_{2},\dots,\omega_{M}\}\subset\Omega$, and the corresponding weights are
$\{w_1,w_2,...,w_{M}\}$ with $\sum_{i=1}^{M}w_i=1$. If the permeability field is given as $\kappa(x,\omega)$,
then we can sample the space of $\Omega$ to obtain realizations.
There are many methods to sample the space $\Omega$, such as Monte Carlo sampling, quasi Monte Carlo sampling \cite{Caflisch:98}, and Latin hypercube sampling. The measure defined on $\Omega_{d}$ is a weighted discrete measure which accurately approximates the
measure on the $\Omega$, that is,
\[
\int_{\Omega_{d}}u(x,\omega)d\omega=\sum_{i}^{M}u(x,\omega_{i})w_{i}\approx\int_{\Omega}u(x,\omega)d\omega.
\]
The fine-grid reference solution $u_{h}\in V_{h}$ is defined by solving
\begin{equation}
\label{eq:fine}
\int_{\Omega_{d}}a(\omega;u,v)d\omega=\int_{\Omega_{d}}l(\omega;v)d\omega,\;\forall v\in V_{h}
\end{equation}
 where $V_{h}$ is the fine-grid finite element space defined by
\[
V_{h}=\{v(x,\omega)\in L^{2}(\Omega_{d};H^{1}(D))\;:\; v(x,\omega_{i})|_{\tau}\in Q_{1}(\tau), \;\forall \tau\in\mathcal{T}^{h},\;\forall\omega_{i}\in\Omega_{d}\},
\]
where $Q_{1}(\tau)$ is a linear function on the element $\tau$. We remark that the fine-grid partition $\mathcal{T}^h$ of the domain $D$ allows us to resolve all heterogeneities of $\kappa(x,\omega)$ in $D$.
In general, solving Eq.(\ref{eq:fine}) is expensive as the dimension of the fine-grid finite element space $V_h$ is large and one needs to solve the multiscale problem for many realizations of $\kappa(x,\omega)$.

In the following sections, we shall design an effective approach that allows one to solve (\ref{eq:fine}) approximately and accurately with much lower computational costs. We shall build a lower dimensional space $V_H$ using the idea of GMsFEM. In order to build the most efficient representation of the solution, we shall not separate the physical variable $x$ and the random variable $\omega$ in our framework. That is, we shall not use a tensor product approach. Indeed, we shall couple the influence in terms of both $x$ and $\omega$ simultaneously in our construction.

\section{Coarse grids}

In this section, we shall introduce the coarse grid discretization of the spacial
domain $D$ and the parameter space $\Omega_{d}$.
One obvious choice of the partition of the product space $D\times \Omega_d$ is to construct the
partitions of $D$ and $\Omega_d$ separately and then use the tensor products.
However, this approach may not be good in practice.
We shall use a more adaptive approach, namely, the partition of the parameter space
depends on the partition of the spatial domain $D$.
In particular, We shall first define a partition for the domain $D$.
Then for each of the element in the partition of $D$, We shall construct a corresponding partition
for the set $\Omega_d$.  Our numerical experiments demonstrate that this new approach is very effective
in solving the heterogenous multiscale problems.

\subsection{Coarse spatial grid}

\begin{figure}[H]
\centering
\includegraphics[scale=0.4]{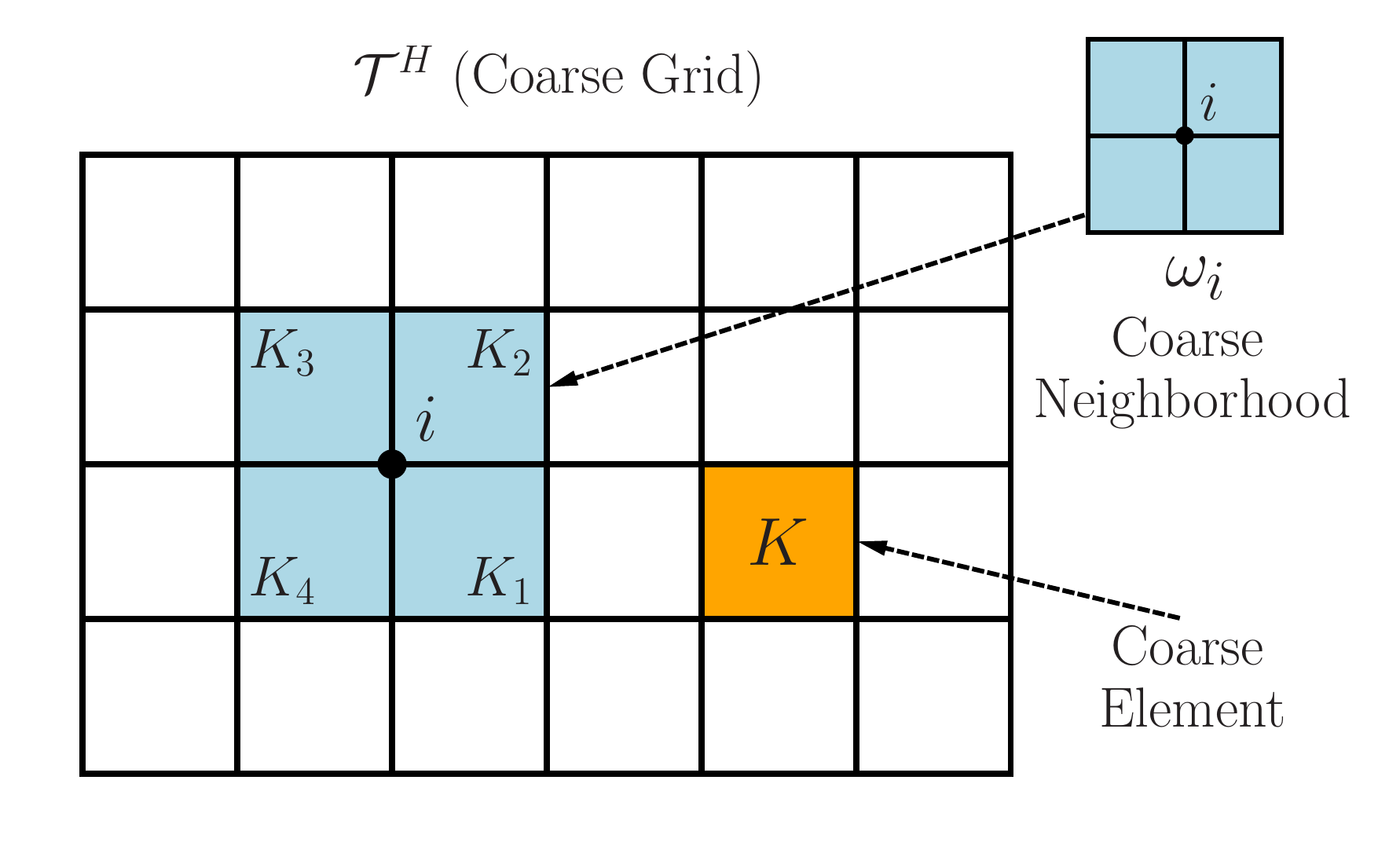}
\protect\caption{Illustration of the coarse grid and related notations.}
\label{illustration}
\end{figure}

We first present the construction of the coarse grid and related notations
for the spatial domain $D$.
We use $\mathcal{T}^{H}$
to denote a coarse grid partition of the domain $D$ and the elements
of $\mathcal{T}^{H}$ are called coarse elements, where $H$ is the
mesh size of the coarse element. To simplify the discussions,
we shall assume that the fine grid partition $\mathcal{T}^h$ is a refinement
of the coarse grid partition $\mathcal{T}^{H}$. Let $x_{i}$, $1\leq i\leq N$, be the
set of interior nodes in the coarse grid, and $N$ is the number of the interior coarse grid nodes. For each coarse node $x_{i}$, the
coarse neighborhood $D_i$ is defined by
\[
D_{i}=\cup\{K_{j}\in\mathcal{T}^{H}:\; x_{i}\in\overline{K}_{j}\},
\]
that is, the union of all coarse elements $K_{j}\in\mathcal{T}^{H}$
having the vertex $x_{i}$. We present an illustration of the course
grid notations in Figure \ref{illustration}.

\subsection{Coarsening $\Omega_d$ via clustering  algorithm.}

Next, we discuss the coarsening of $\Omega_d$ via clustering algorithm. The clustering  algorithm will be performed for each neighborhood of $D_i$ in general. In simplified cases, one can assume the same clustering for all $D_i$'s; however, we emphasize that in general cases, one needs to coarsen $\Omega_d$ for each coarse spatial grid block $D_i$ separately.

To perform the clustering, we shall use a distance function $d^i(\cdot,\cdot)$ defined on $\Omega_d$ for the coarse neighborhood $D_i$. The distance function $d^i(\cdot,\cdot)$ is defined in the solution space instead of the sample space. More precisely, for any two elements $\omega_n,\omega_m \in \Omega_d$,
we define the distance $d^i(\omega_n,\omega_m)$
be the distance of the two functions $\phi_n,\phi_m$ in $D_i$,
where $\phi_n,\phi_m$ are solutions of $-\nabla\cdot(\kappa(x,\omega)\nabla \phi(x,\omega))=0$
with $\omega=\omega_n, \omega_m$ in $D_i$ respectively and with an appropriate boundary condition on $\partial D_i$. However, solving the multiscale problem to get $\phi_n$ and $\phi_m$ for every pair of elements $\omega_n$ and $\omega_m$ in $\Omega_d$ is expensive. Therefore, we present a simplified way to reduce the computational cost in our clustering algorithm.

\begin{figure}[H]
\centering
\includegraphics[scale=0.4]{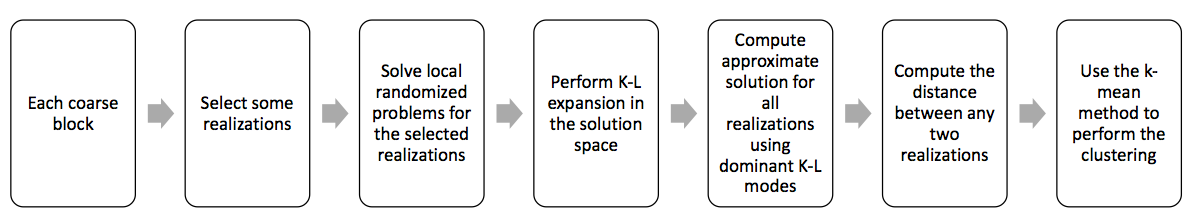}
\protect\caption{Illustration of the clustering algorithm}
\label{illustration1}
\end{figure}

For each coarse neighborhood $D_{i}$, we first construct a snapshot space for a subset in $\Omega_d$. We denote this subset of $\Omega_d$ by $\Omega_d^{\text{subset}}$. Our objective is to find a simplified distance
functions in $\Omega_d$ that can be used for clustering. The main idea is to use randomized harmonic extensions
and KL expansion. The clustering algorithm is illustrated in Figure \ref{illustration1}.

We consider the $i$-\emph{th} coarse neighborhood $D_i$ and a fixed subset $\Omega_d^{\text{subset}}$ of $\Omega_d$. We let $D_i^+$ be an oversampled domain for $D_i$, that is, $D_i \subset D_i^+$. More details about the oversampling technique for multiscale finite element method can be found in \cite{HouEfendievWu:2000} and references therein. In addition, we let $V_h^{i,+}$ be the subspace of $V_h$ defined by
\[
V_{h}^{i,+}=\{v(x,\omega)\in L^{2}(\Omega_{d}^{\text{subset}};H^{1}(D_i^+))\;:\; v(x,\omega)|_{\tau}\in Q_{1}(\tau), \;\forall \tau\in\mathcal{T}^{h}\cap D_i^+,\;\forall\omega\in\Omega_{d}^{\text{subset}}\},
\]
and the subspace $V_{h,0}^{i,+}$ be the subspace of $V_h^{i,+}$ containing functions with zero trace on $\partial D_i^+$.
Then we construct
a set of snapshot functions $\psi_{j}^{i}(x)$, by solving the local problem in the oversampling sub-domain
$D_{i}^{+}\times \Omega_d^{\text{subset}}$. In particular, we find $\psi_j^i(x) \in V_h^{i,+}$ such that
\begin{equation}
\label{eq:CellProblemRandomBC}
a(\omega;\psi_{j}^{i}(x,\omega),v)=l(\omega,v), \quad \forall v\in V_{h,0}^{i,+},  \quad \forall\omega\in\Omega_{d}^{\text{subset}},
\end{equation}
subject to the boundary conditions $\psi_{j}^{i}(x)|_{\partial D_{i}^{+}}=R_{j}$, $j=1,...,k^{i}$,
where $R_{j}$ is a discrete random function defined with respect to the fine grid boundary
point on $\partial D_{i}^{+}$ and $k^{i}$ are the number of local problem that one needs to solve
in coarse neighborhood $D_i$.

How to choose the appropriate number of snapshot functions $\psi_{j}^{i}(x)$ number $k^{i}$  will effect the performance of our method. Our goal is to construct a few effective boundary conditions so that the snapshot functions $\psi_{j}^{i}(x)$ can be used to approximate the solution space. We find that this is closely related to the range-finding problem. In discrete level, we can formulate the range-finding problem as: for a matrix $T$, we want to find a matrix $Q$ with independent columns such that its column space
accurately approximates the column space $T$, i.e.,
\begin{equation}
\label{eq:ApproximateOperator}
||T-QQ^{T}T||\leq \epsilon,
\end{equation}
where $\epsilon$ is a small parameter and $||\cdot||$ is a matrix norm.

Randomized range finding algorithm is very efficient to construct the matrix $Q$  and to capture the column space of $T$. See \cite{Tropp:11,Rokhlin:07} for more details about the randomized algorithms for finding low-dimensional structures. To justify the approximation error $||T-QQ^{T}T||$, we denote $B=(I-QQ^{T})T$ and
draw $k$ random vectors $\omega^{j}$, $j=1,...,k$. We use $T$ to act on these $\omega^{j}$, $i=1,...,k$, and compute the residual after projecting the images  $T\omega^{j}$ to the column space of $Q$, namely $||(I-QQ^{T})T\omega^{j}||$. If all the $k$ number of residuals are smaller than $\frac{\epsilon}{\alpha}\sqrt{\frac{\pi}{2}}$, then the condition \eqref{eq:ApproximateOperator} holds except with a small probability $\alpha^{-k}$, where $\alpha>1$.

The good thing about the randomized range-finding algorithms and the error estimation \eqref{eq:ApproximateOperator} is that they do not require access to each entry of the matrix $T$, but only require the matrix-vector multiplication. For our problem, this matrix-vector multiplication corresponding to solve \eqref{eq:CellProblemRandomBC} with certain  random boundary conditions. Therefore, we employ the randomized range finding philosophy and the error estimation in our problem to determine an optimal number $k^{i}$ so that we can reduce the computational cost.


In the next step, we perform a KL expansion of the snapshot functions $\{ \psi_j^i(x,\omega)\}$ and
obtain the dominant modes in the solution space
\begin{equation}
\label{eq:kle_snapshotfunctions}
\psi_{j}^{i}(x,\omega)=\bar{\psi}_{j}^{i}(x)+\sum_{l \geq 1}p_{j,l}^{i}(\omega)\phi_{l}^{i}(x), 
\end{equation}
where $\bar{\psi}_{j}^{i}(x)=\cfrac{1}{|\Omega_d^{\text{subset}}|}\int_{\Omega_d^{\text{subset}}} \psi_{j}^{i}(x,\omega) \; d\omega$ is the mean of snapshot functions and $|\Omega_d^{\text{subset}}|$ is the number of sample in $\Omega_d^{\text{subset}}$.  We note that
$\phi_{l}^{i}(x)\in V_{h,0}(D_{i}^{+})$ are defined by
\[
V_{h,0}(D_i^+) =\{v(x)\in H^{1}_0(D_i^+)\;:\; v(x)|_{\tau}\in Q_{1}(\tau), \;\forall \tau\in\mathcal{T}^{h}\cap D_i^+ \}.
\]
We also note that $\bar{\psi}_{j}^{i}(x)$ satisfies the same boundary condition as $\psi_{j}^{i}(x)$.

Among the set $\{ \phi_l^i(x), l=1,2,...\}$ of the KL expansion of the snapshot functions in \eqref{eq:kle_snapshotfunctions}, we shall restrict them in the coarse neighborhood $D_i$ and then
take the first $L^i$ dominant parts and form the reduced spatial snapshot space $V_{\text{snap}}^i$.
The low-dimensional reduced space $V_{\text{snap}}^i$ is defined by
\[
V_{\text{snap}}^i = \text{span}\{\phi_{l}^{i}(x),\; 1\leq l\leq L^{i}\}.
\]
The key idea of our new algorithm is to use this low-dimensional space to solve the local problem for each realization in $\Omega_d$.
More precisely, for each $\omega \in \Omega_d$,
we find a function $\tilde{\psi}_{j}^{i}(x, \omega)$ such that
\begin{equation*}
\tilde{\psi}_{j}^{i}(x,\omega)=\bar{\psi}_{j}^{i}(x)+\sum_{1\leq l\leq L^{i}}\tilde{p}_{j,l}^{i}(\omega)\phi_{l}^{i}(x), \quad\quad x\in D_i,
\end{equation*}
and
\[
a(\omega,\tilde{\psi}_{j}^{i},v)=l(\omega,v), \quad\quad \forall v\in V_{\text{snap}}^{i}.
\]
We remark that solving the above problems is very efficient due to the small dimension of the solution space $V_{\text{snap}}^i$,
which enables us  to construct a distance function in $\Omega_d$. To be specific, we use local solutions to define a distance function on $\Omega_d$, denoted as $d^{i}(\cdot,\cdot):\Omega_d\times\Omega_d\rightarrow\mathbb{R}^{+}\bigcup\{0\}$, by
\begin{equation}
d^{i}(\omega_{n},\omega_{m})=\sqrt{\sum_{j}\sum_{1\leq l\leq L^{i}}(\tilde{p}_{j,l}^{i}(\omega_{n})-\tilde{p}_{j,l}^{i}(\omega_{m}))^2}. \label{eq:DistanceFunc_ReducedSpace}
\end{equation}
We then use the k-means clustering algorithm \cite{Hartigan:1979} with the distance function defined in \eqref{eq:DistanceFunc_ReducedSpace} to cluster
the sampling space $\Omega_d$ into $J^i$ clusters $\Omega^{i}_{j}$, $j=1,...,J^i$,
such that, $\Omega_d=\cup_{1\leq j\leq J^i}\Omega_{j}^{i}$. 

\begin{remark}
We should point out that the essential idea is that the clustering of the
sampling space $\Omega_d=\cup_{1\leq j\leq J^i}\Omega_{j}^{i}$ depends on the spatial location $D_i$. This
enables us to efficiently explore the heterogeneities of $\kappa(x,\omega)$ as well as the solution space.
\end{remark}

\section{The construction of offline space}

In this section, we describe the construction of the local offline basis functions. Let $D_{i}$, $i=1,...,N$ and $\Omega_{j}^{i}$, $j=1,...,J^i$ be a given coarse neighborhood and a cluster of $\Omega_d$, respectively. In the constructing  process, we first construct a snapshot space $V_{\text{snap}}^{(i,j)}$ for $D_{i}\times\Omega_{j}^i$. In the notation $V_{\text{snap}}^{(i,j)}$, the indices $i$ are corresponding to coarse neighborhood $D_i$ and the indices $j$ are induced by $i$ when we use the clustering algorithm.

The construction of the snapshot space $V_{\text{snap}}^{(i,j)}$ consists of solving the local problems for several random boundary conditions. To construct the offline space $V_{\text{off}}^{(i,j)}$, we need to
solve a local spectral problem in the $V_{\text{snap}}^{(i,j)}$ for
the dimension reduction. We shall discuss the full description of the
construction in the following subsections.

\subsection{Construct the snapshot space}

The definition of the snapshot space $V_{\text{snap}}^{(i,j)}$ is based on $a(\omega)$-harmonic extensions. Let $B_{h}(D_{i})$ be the set of all fine-grid nodes lying on the boundary of the coarse neighborhood of $\partial D_{i}$. For each
fine grid node $x_{i}$, the discrete delta function $\delta_{k}^{h}(x)$
is defined by
\[
\delta_{k}^{h}(x_{l})=\begin{cases}
1, & k=l\\
0 & k\neq l
\end{cases}\quad,x_{l}\in B_{h}(D_{i}).
\]

Next, we define the $k$-th snapshot basis function $\psi_{k}^{(i,j),\text{snap}}\in V_{h}(D_{i})$
as the solution of
\begin{align}
\int_{D_{i}}\bar{\kappa}^{(i,j)}(x)\nabla\psi_{k}^{(i,j),\text{snap}}\cdot\nabla v & =0,\;\forall v\in V_{h,0}(D_{i})
\label{eq:k-thSnapshotBasisFunction_Eq} \\
\psi_{k}^{(i,j),\text{snap}}(x) & =\delta_{k}^{h}(x),\;\text{on }\partial D_{i},
\label{eq:k-thSnapshotBasisFunction_BC}
\end{align}
where $\bar{\kappa}^{(i,j)}(x) = (\cfrac{1}{|\Omega^{i}_{j}|}\int_{\Omega^{i}_{j}}\kappa(x,\omega)d\omega)$.
Here, we take the mean of $\kappa(x,\omega)$ within a cluster for computing offline spaces. In general, one can use a typical realization or a few realizations in a similar fashion. The main idea is that within each cluster, we assume that the realizations are similar in an appropriate metric (i.e., the distance between
randomized snapshots).

The dimension of $V_{\text{snap}}^{(i,j)}$ is equal to the number
of fine grid nodes lying on $\partial D_{i}$. We can use the randomized snapshots with oversampling technique to reduce the dimension of $V_{\text{snap}}^{(i,j)}$ and therefore reduce the computational cost of solving these snapshot basis functions.

\subsection{Construct the offline space}

To obtain the offline space $V_{\text{off}}^{(i,j)}$, we shall perform a spectral decomposition in
the snapshot space.  The spectral decomposition is defined as to find $(\phi^{(i,j)(x)}_k,\lambda^{(i,j)}_k ) \in V_{\text{snap}}^{(i,j)}\times\mathbb{R}$ such that
\begin{equation}
\int_{D_{i}}\bar{\kappa}^{(i,j)}(x)\nabla\phi_{k}^{(i,j)}(x)\cdot\nabla v  =
\lambda^{i,j}_k \int_{D_{i}}\bar{\kappa}^{(i,j)}(x)|\nabla \chi_i(x)|^2\phi_{k}^{(i,j)}(x) v ,\;
\forall v\in V_{\text{snap}}^{(i,j)},
\label{eq:eigen_decomposition}
\end{equation}
where $\{\chi_i(x)\}$ is the partition of unity functions for the spatial domain $D$ with respect to the partition $\{ D_i \}$
and $\bar{\kappa}^{(i,j)}(x)$ is defined the same as in \eqref{eq:k-thSnapshotBasisFunction_Eq}. We assume the eigenvalues are arranged in an ascending order. Then, the offline basis functions for the coarse neighborhood $D_i\times \Omega^{i}_{j}$ are defined by $\tilde{\phi}^{(i,j)}_{k}(x,\omega) = \chi_i(x) I_{\Omega^{i}_{j}}(\omega)\phi^{(i,j)}_{k}(x)$, for $1\leq k \leq M^{(i,j)}$, where $I_{\Omega^{i}_{j}}$ is the characteristic function of $\Omega^{i}_{j}$. That is,
\begin{equation}
V_{\text{off}}=\text{span}\{ \tilde{\phi}^{(i,j)}_{k}: \forall i, 1\leq j \leq J^i, 1\leq k\leq M^{(i,j)}\}.
\label{eq:OfflineBasisFunctions}
\end{equation}

\subsection{Global formulation}
The offline space $V_{\text{off}}$ can be used to solve the multiscale stochastic problem (\ref{eq:PDE})
for any input parameter $\omega$ and any right hand side function $f(x,\omega)$. In particular,
the stochastic GMsFEM for the multiscale stochastic problem (\ref{eq:PDE}) is to find $u_{\text{ms}} \in V_{\text{off}}$ such that
\begin{equation}
\label{eq:scheme}
\int_{\Omega_d} a(\omega;u_{\text{ms}},v) \; d\omega=\int_{\Omega_d}l(\omega;v) \; d\omega, \quad\quad \forall v\in V_{\text{off}}.
\end{equation}
By solving \eqref{eq:scheme} using Galerkin method, we can obtain the numerical approximation solution $u_{\text{ms}}$ to the multiscale stochastic solution.

\subsection{Online basis functions}\label{sec:online}

The numerical solution obtained in the offline space has already produced a good approximation to the multiscale stochastic solution. In some cases, however, we need to improve the accuracy of the solution by using some online basis functions,
which are constructed in the online stage.

To demonstrate the main idea of the online update, we consider a coarse neighborhood $D_i$ and the set $\Omega_j^i$ for example. We define the following residual function
\begin{equation}
R_{i,j}(\omega; v) = \int_{\Omega_j^i} a(\omega;u_{\text{ms}},v) \; d\omega - \int_{\Omega_j^i}l(\omega;v) \; d\omega, \quad \forall v \in L^2(\Omega_j^i; V_{h,0}(D_i)),
\end{equation}
where
\[
V_{h,0}(D_i) =\{v(x)\in H^{1}_0(D_i)\;:\; v(x)|_{\tau}\in Q_{1}(\tau), \;\forall \tau\in\mathcal{T}^{h}\cap D_i \}.
\]
To construct an online basis function $\phi^{(i,j)}_{\text{on}} \in L^2(\Omega_j^i; V_{h,0}(D_i))$, we solve the following problem
\begin{equation}
\label{eq:online}
\int_{\Omega_j^i} a(\omega;\phi^{(i,j)}_{\text{on}},v) \; d\omega =
R_{i,j}(\omega; v), \quad \forall v \in L^2(\Omega_j^i; V_{h,0}(D_i)).
\end{equation}
We assume a piecewise constant approximation in the space $\Omega_j^i$. The above problem (\ref{eq:online}) leads to
\begin{equation}
\label{eq:online1}
a(\omega_k;\phi^{(i,j)}_{\text{on}}(k,\cdot),v)  =
R_{i,j}(\omega_k; v), \quad \forall v \in V_{h,0}(D_i).
\end{equation}
The basis function $\phi^{(i,j)}_{\text{on}}$ is added to the space $V_{\text{off}}$. In general, we do not need to solve (\ref{eq:online1}) for every $\omega_k$ in the set $\Omega_j^i$. We only need to choose those $\omega_k$ such that the residual $R_{i,j}(\omega_k;v)$ are large and update its corresponding reduced basis space.

\subsection{Discussions}
Before we end this section, we give some general discussions about our algorithms as follows.
\begin{itemize}
\item In our algorithm, we have used $k$-mean algorithm for clustering and
a multiscale approximation for computing distances for all pairs. In general,
we can use other clustering algorithms \cite{manning2008introduction}
and more accurate approximation for distances. This can make the coarsening of the uncertainty space (clustering) to provide a better convergence and we plan to study this in our future work.

\item In the paper, we do not analyze the accuracy of the method. We expect
the accuracy depends on the realization-based multiscale approximations  and the cluster sizes. The first one has been analyzed in our previous papers \cite{chung2014adaptive} and can be improved \cite{chung2017constraint}. The error due to the cluster size will be in the form of the variance of the solution space, which we plan to study in our future works.

\item The proposed methods can take an advantage of the adaptivity in spatial space
and uncertainty space. In particular, the number of multiscale
basis functions in each coarse-grid block can be controlled via some error
indicators \cite{chung2014adaptive}. Adaptivity for clustering will involve changing the cluster
sizes and adding more or less elements in each cluster. This can be expensive in general.

\item The online step is important to accelerate the convergence. By adding the
online basis functions adaptively, we can decrease the error significantly at a cost
of reformulating the coarse-grid problem. We believe the accuracy of the online
method depends on the offline spaces as in \cite{chung2015residual}.

\end{itemize}

\section{Numerical Results}
In this section, we perform some numerical experiments to test the performance and accuracy of the proposed
cluster-based generalized multiscale finite element method for multiscale elliptic PDEs with random coefficients. As we will demonstrate in our paper, our new method could offer accurate numerical solutions with considerable computational savings over traditional stochastic methods.
In the results below, we will use two types of error quantities to show the performance of our method. The
first type of errors is
\begin{equation*}
e_{1,\Omega} := \sqrt{\int_{\Omega}\int_{D}|u_{h}-u_{H}|^{2}dxd\omega},\quad
e_{2,\Omega} := \sqrt{\int_D |\int_{\Omega} (u_{h}-u_{H}) d\omega|^2 dx}
\end{equation*}
which measures the $L^2$ norm of the error over all realizations. The second type of errors is
\begin{equation*}
e_{1,S} := \sqrt{\sum_{i=1}^{M}\int_{D}|u_{h}(x,\omega_{i})-u_{H}(x,\omega_{i})|^{2}dxd\omega}, \quad
e_{2,S} := \sqrt{\int_D |\sum_{i=1}^{M} u_{h}(x,\omega_{i})-u_{H}(x,\omega_{i}) d\omega|^{2}dx}
\end{equation*}
which measures the $L^2$ norm of the error over $M$ realizations, where $M=10$ in our numerical experiments.

\subsection{Permeability with high-contrast inclusions and channels} \label{sec:high-contrast}
\begin{figure}[ht]
\centering
\includegraphics[scale=0.48]{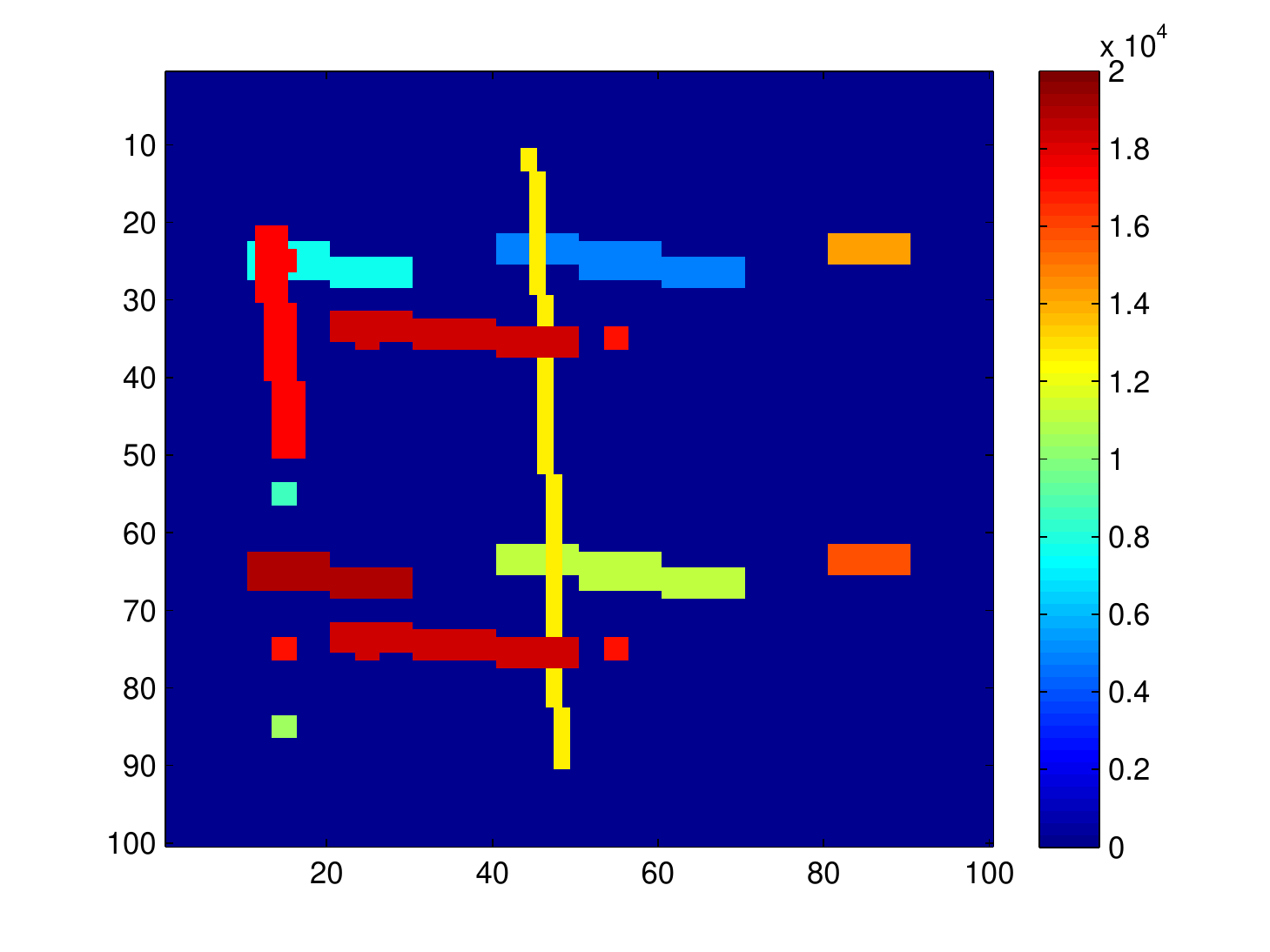}\includegraphics[scale=0.48]{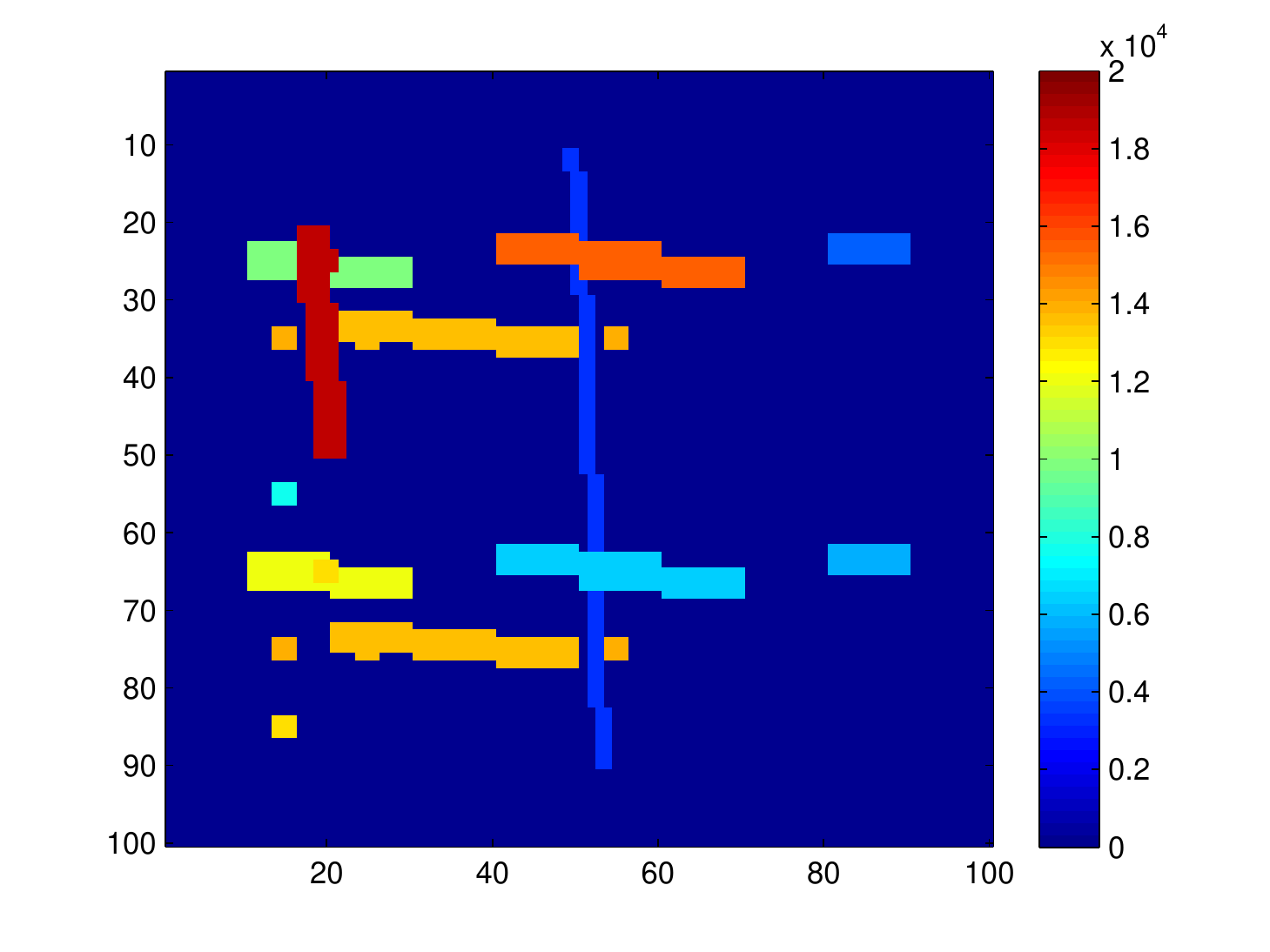}
\caption{Permeability coefficient $\kappa$ of two different realizations for case 1.}
\label{fig:medium_case1}
\end{figure}

In this example, we consider a permeability containing high-contrast inclusions and channels (see Figure \ref{fig:medium_case1}). More precisely, we consider the following 2D elliptic PDE with
random coefficient,
\begin{align}
\label{eq:numericalexample1}
-\nabla\cdot(a(x,\omega)\nabla u(x,\omega))&=f(x),  \quad  (x,\omega)\in D\times\Omega  \\
u(x,\omega)&=0, \quad  (x,\omega)\in \partial D\times\Omega
\end{align}
To introduce randomness, we move these high-contrast inclusions at different directions
and change their permeability values, while keeping the values at a high
range. In Figure \ref{fig:medium_case1}, we depict two realizations.

In Table \ref{tab:case1_cluster1}, we show numerical results for one of the clusters. In this table, we compare $L^2$ norms
of the mean and realizations between the fine-grid solution and the coarse-grid solution across the space
and realizations. In particular, we show the results when choosing all realizations and when selecting $10$ random realizations. In our results, we increase the number of basis functions and add more basis functions in each
coarse-grid block. The latter can be done adaptively. We observe that the error has been reduced to $8$\% when $5$ basis functions are used, which is similar to the case of using $3$ basis functions. We note that this is the error due to the clustering size and can not be reduced further unless we change the cluster size.

In the next tables, we increase the number of cluster. In Table \ref{tab:case1_cluster3}, we use $3$ clusters and increase
the number of basis functions per coarse element. We observe from this table that the error decreases as we increase the number of basis functions and, moreover, the error is smaller when using $3$ clusters compared to the case of using $1$ cluster.

Next, in Table \ref{tab:case1_cluster5}, we use $5$ clusters and increase the number of basis functions per coarse element.
In this case, we again observe an error decreasing as the number of basis is being increased and smaller error compared to the case with $3$ clusters. In this table, we also show an error corresponding when using GMsFEM basis functions for each realization. The latter is used to identify the error due to GMsFEM coarse-grid discretization. We observe that this error is about $4$\% compared to $6$\%, which indicates that our approach requires a few clusters in this example to achieve a good accuracy.

\begin{table}[ht]
\centering
\begin{tabular}{|c|c|c|}
\hline
number of basis  & $e_{1,\Omega}$ & $e_{2,\Omega}$\tabularnewline
\hline
1 & 24.38\% & 24.19\%\tabularnewline
\hline
3 & 14.57\% & 14.22\%\tabularnewline
\hline
5 & 12.28\% & 11.86\%\tabularnewline
\hline
\end{tabular}

\centering
\begin{tabular}{|c|c|c|}
\hline
number of basis  & $e_{1,S}$ & $e_{2,S}$\tabularnewline
\hline
1 & 21.76\% & 21.59\%\tabularnewline
\hline
3 & 10.89\% & 10.57\%\tabularnewline
\hline
5 & 8.73\% & 8.39\%\tabularnewline
\hline
\end{tabular}

\label{tab:case1_cluster1}
\caption{Errors for the case $1$ using \#cluster = 1.}
\end{table}

\begin{table}[ht]
\centering
\begin{tabular}{|c|c|c|}
\hline
number of basis  & $e_{1,\Omega}$ & $e_{2,\Omega}$\tabularnewline
\hline
1 & 23.96\% & 23.76\%\tabularnewline
\hline
3 & 12.50\% & 12.06\%\tabularnewline
\hline
5 & 10.25\% & 9.75\%\tabularnewline
\hline
\end{tabular}

\begin{tabular}{|c|c|c|}
\hline
number of basis  & $e_{1,S}$ & $e_{2,S}$\tabularnewline
\hline
1 & 21.76\% & 21.59\%\tabularnewline
\hline
3 & 9.52\% & 9.15\%\tabularnewline
\hline
5 & 7.37\% & 7.01\%\tabularnewline
\hline
\end{tabular}
\caption{Errors for the case 1 using \#cluster = 3.}
\label{tab:case1_cluster3}
\end{table}

\begin{table}[ht]
\begin{centering}
\begin{tabular}{|c|c|c|}
\hline
number of basis  & $e_{1,\Omega}$ & $e_{2,\Omega}$\tabularnewline
\hline
1 & 23.90\% & 23.70\%\tabularnewline
\hline
3 & 12.26\% & 11.86\%\tabularnewline
\hline
5 & 8.82\% & 8.35\%\tabularnewline
\hline
\end{tabular}
\par\end{centering}

\begin{centering}
\begin{tabular}{|c|c|c|}
\hline
number of basis  & $e_{1,S}$ & $e_{2,S}$\tabularnewline
\hline
1 & 21.76\% & 21.59\%\tabularnewline
\hline
3 & 9.59\%(5.64\%) & 9.24\%(5.58\%)\tabularnewline
\hline
5 & 6.31\%(3.92\%) & 6.10\%(3.88\%)\tabularnewline
\hline
\end{tabular}
\par\end{centering}
\centering{}\protect\caption{Errors for the case 1 using \#cluster = 5. The error in the parenthesis
shows the error using GMsFEM basis function using the specific realization.}
\label{tab:case1_cluster5}
\end{table}

\begin{figure}[ht]
\centering
\includegraphics[scale=0.6]{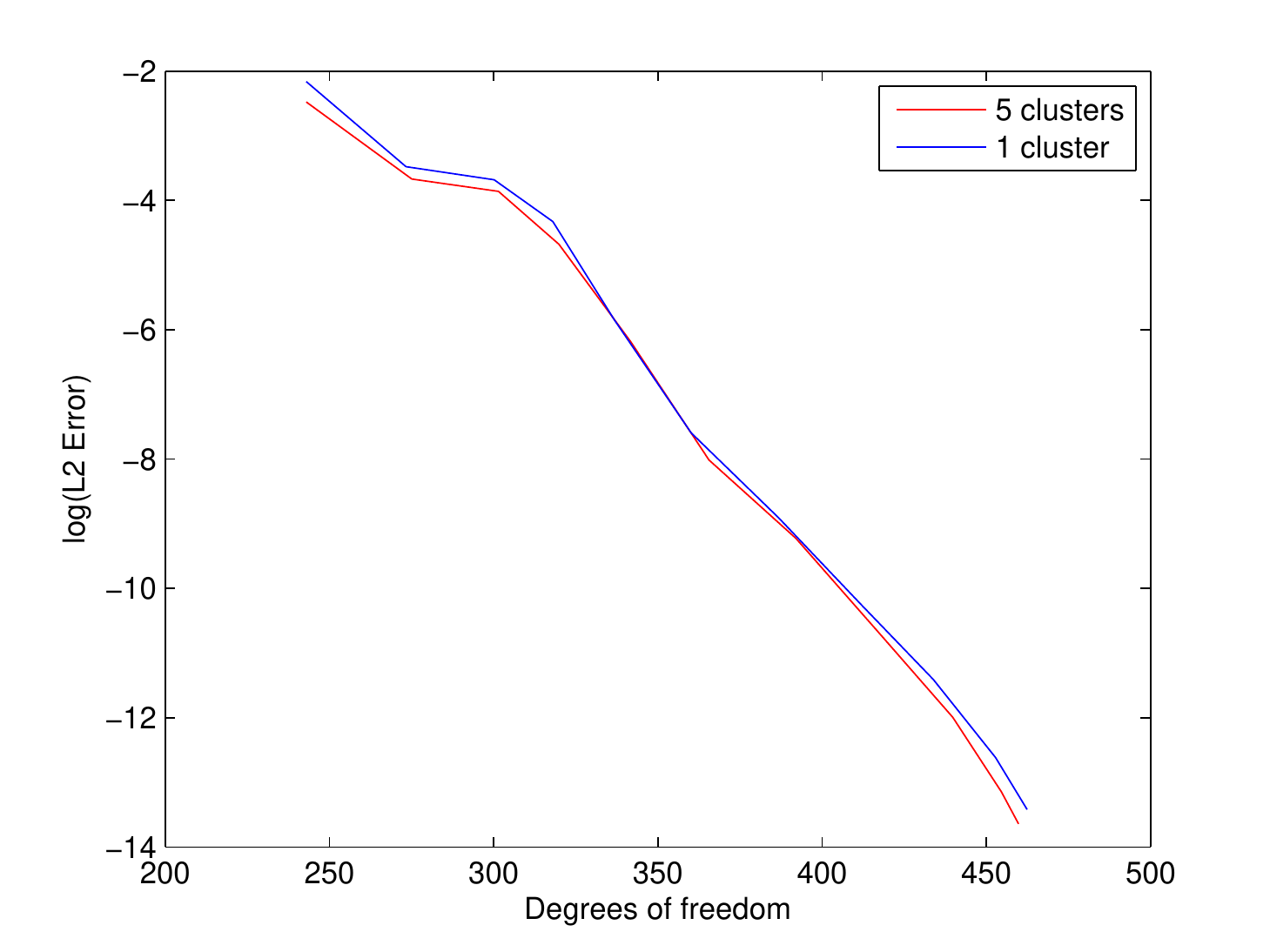}
\caption{Error decay for online basis functions using 3 offline basis functions with different numbers of clusters.}
\label{fig:medium_case1_OnlineResult}
\end{figure}

Finally, we present a test case to show the performance of using online basis functions defined in Section \ref{sec:online}.
In Table \ref{tab:online}, we show the numerical results.
In our computations, we start the process by using $3$ offline basis functions, and this step corresponds to zeroth
online iteration. Then we show the errors for the next $3$ online iterations, where one online basis is added
for each iteration. We observe that the method is able to produce very accurate results after a couple of online iterations.
In Figure~\ref{fig:medium_case1_OnlineResult}, we present the convergence history the
online iteration, namely, we show the logarithm of the error versus the number of degrees of freedoms.
We observe clearly the exponential decay of the error from this figure.

\begin{table}[ht]
\centering
\begin{tabular}{|c|c|c|}
\hline
\# online iteration  & $e_{1,S}$ & $e_{2,S}$\tabularnewline
\hline
0 & 11.53\% & 11.31\%\tabularnewline
\hline
1 & 3.08\% & 2.55\%\tabularnewline
\hline
2 & 2.52\% & 1.96\%\tabularnewline
\hline
3 & 1.32\% & 0.89\%\tabularnewline
\hline
\end{tabular}

\begin{tabular}{|c|c|c|}
\hline
\# online iteration  & $e_{1,S}$ & $e_{2,S}$\tabularnewline
\hline
0 & 8.40\% & 8.20\%\tabularnewline
\hline
1 & 2.55\% & 1.94\%\tabularnewline
\hline
2 & 2.11\% & 1.48\%\tabularnewline
\hline
3 & 0.93\% & 0.48\%\tabularnewline
\hline
\end{tabular}
\caption{Online results for using 3 offline basis functions with different numbers of clusters. Top: results for 1 cluster. Bottom: results for  5 clusters.}
\label{tab:online}
\end{table}

\subsection{Permeability with randomness and spatial heterogeneities} \label{sec:heterogeneities}

\begin{figure}[ht]
\centering
\includegraphics[scale=0.48]{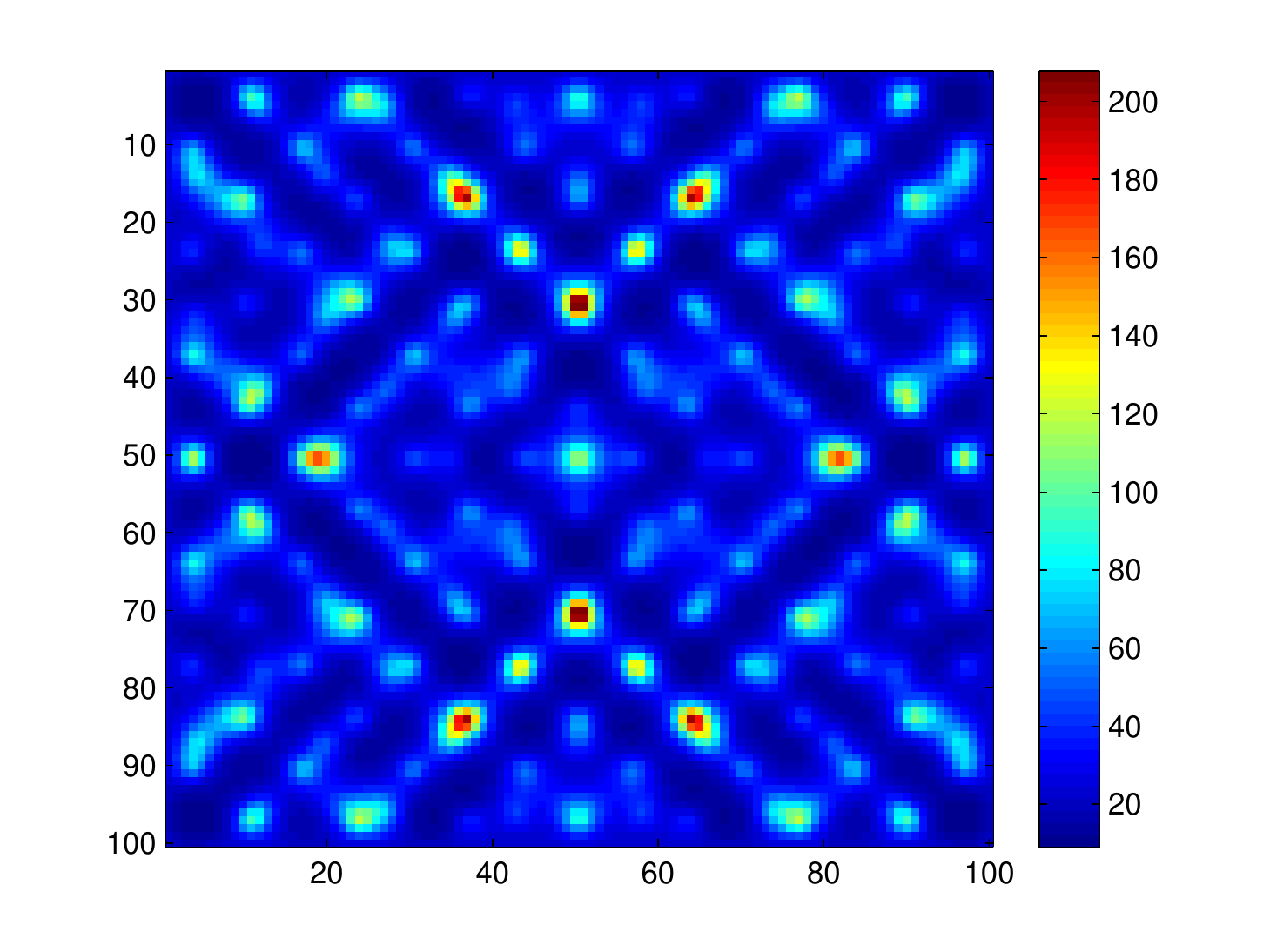}\includegraphics[scale=0.48]{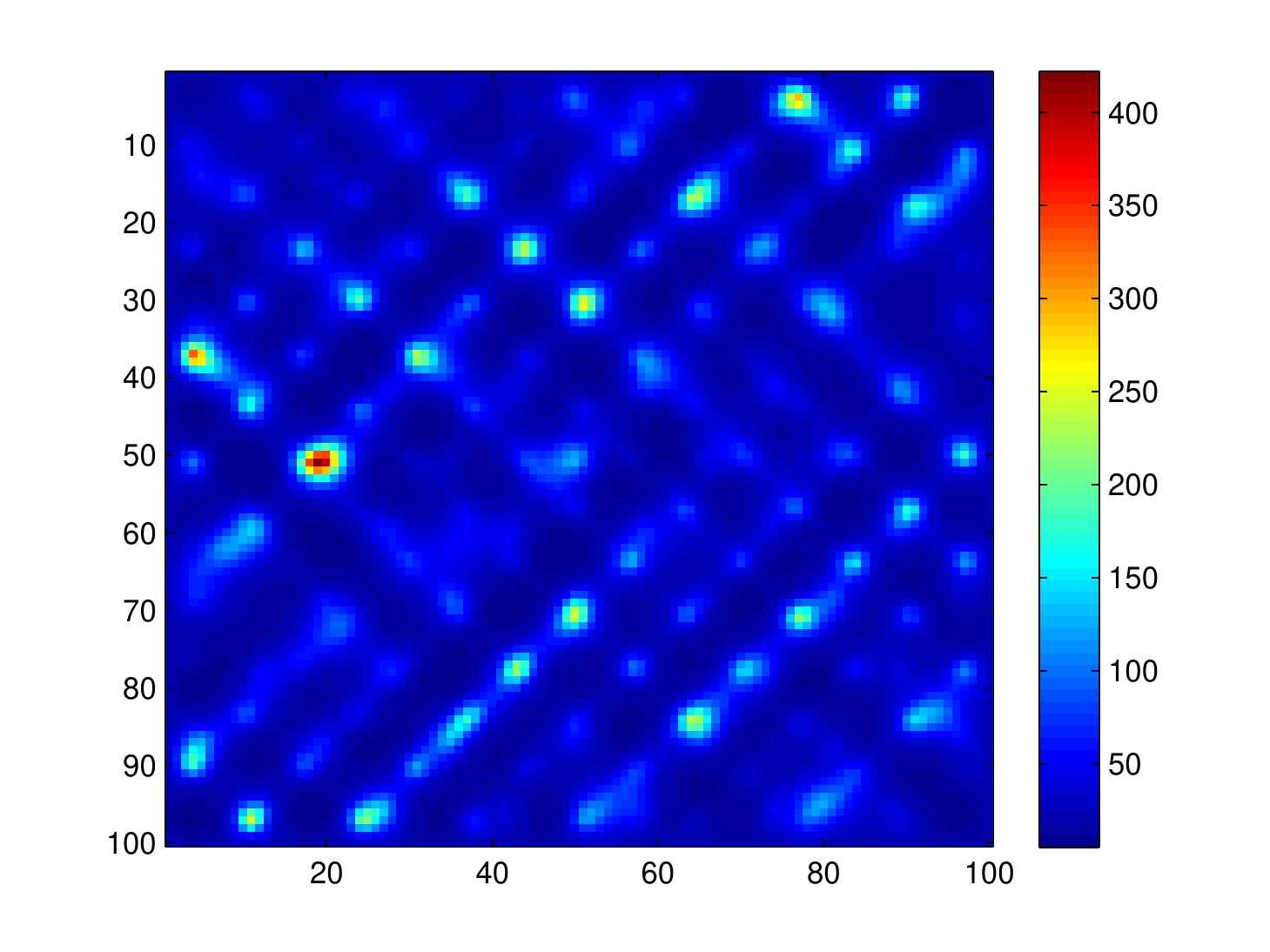}
\caption{Permeability coefficient $\kappa$ of two different realizations for Eq. \eqref{eq:CoefficientForEx2}.}
\label{fig:medium_case2}
\end{figure}

In this example, we consider the same 2D elliptic PDE \eqref{eq:numericalexample1}
with homogeneous boundary condition and a random permeability field
\begin{align}
\label{eq:CoefficientForEx2}
a(x,\omega)=e^{0.1+\frac{2+\sin(7\pi x_{1})\sin(8\pi x_{2})}{2+\sin(9\pi x_{1})\sin(7\pi x_{2})}\xi_{1}(\omega)
+\frac{2+\sin(13\pi x_{1})\sin(11\pi x_{2})}{2+\sin(11\pi x_{1})\sin(13\pi x_{2})}\xi_{2}(\omega)
+\frac{2+\sin(12\pi x_{1})\sin(14\pi x_{2})}{2+\sin(15\pi x_{1})\sin(15\pi x_{2})}\xi_{3}(\omega)},
\end{align}
where we have $3$ random variables $\xi_{i}(\omega)$ follow the standard normal distribution $\mathcal{N}(0,1)$ that represent the randomness and spatial heterogeneities. We depict two realizations of this permeability field
in Figure \ref{fig:medium_case2}.

As in our previous example, We shall vary the number of clusters and the number of basis functions.
In Tables \ref{tab:case2_cluster6} and \ref{tab:case2_cluster10}, we consider two cluster sizes with $6$ and $10$ clusters.
In each case, we vary the number of basis functions in each coarse-grid block. We observe from these numerical results that
the method based on $6$ clusters has provided an accurate result that is very comparable to that obtained from using $10$ clusters.

In Table \ref{tab:case2_cluster10}, we present the numerical results when using the multiscale basis functions for each
realization (see the numbers in the parentheses). These results show that $10$ clusters provide almost the same results
compared when using multiscale basis functions for each realization and this, it is sufficient to have only 10 clusters for our multiscale simulations.

\begin{remark}
We should point out that the coefficient $a(x,\omega)$ in \eqref{eq:CoefficientForEx2} is a highly nonlinear functional of $\xi_{i}(\omega)$ and does not have the affine parameter dependence property. The clustering algorithm proposed in this paper helps us automatically explore the heterogeneities of the solution space and constructs the reduced basis functions. Thus, our method can efficiently solve this type of challenging problem.
\end{remark}

\begin{table}[H]
\centering
\begin{tabular}{|c|c|c|}
\hline
number of basis  & $e_{1,\Omega}$ & $e_{2,\Omega}$\tabularnewline
\hline
1 & 17.33\% & 17.26\%\tabularnewline
\hline
3 & 11.24\% & 11.14\%\tabularnewline
\hline
5 & 8.70\% & 8.59\%\tabularnewline
\hline
\end{tabular}

\begin{tabular}{|c|c|c|}
\hline
number of basis  & $e_{1,S}$ & $e_{2,S}$\tabularnewline
\hline
1 & 16.91\% & 16.80\%\tabularnewline
\hline
3 & 10.64\% & 10.54\%\tabularnewline
\hline
5 & 8.20\% & 8.10\%\tabularnewline
\hline
\end{tabular}
\caption{Errors for the case 2 using \#cluster = 6.}
\label{tab:case2_cluster6}
\end{table}

\begin{table}[H]
\centering
\begin{tabular}{|c|c|c|}
\hline
number of basis  & $e_{1,\Omega}$ & $e_{2,\Omega}$\tabularnewline
\hline
1 & 17.26\% & 17.21\%\tabularnewline
\hline
3 & 10.11\% & 10.02\%\tabularnewline
\hline
5 & 8.01\% & 7.92\%\tabularnewline
\hline
\end{tabular}

\begin{tabular}{|c|c|c|}
\hline
number of basis  & $e_{1,S}$ & $e_{2,S}$\tabularnewline
\hline
1 & 16.96\% & 16.86\%\tabularnewline
\hline
3 & 9.29\% & 9.23\%\tabularnewline
\hline
5 & 7.83\%(7.18\%) & 7.80\%(7.13\%)\tabularnewline
\hline
\end{tabular}
\caption{Errors for the case 2 using \#cluster = 10.}
\label{tab:case2_cluster10}
\end{table}

%
%
%

\section{Concluding remarks}

In the paper, we propose a generalized multiscale approach for stochastic problems. Though construction of multiscale basis functions on a coarse spatial grid is well studied, there has not been much research on coarsening uncertainty
space for multiscale problems. The coarsening of the uncertainty space needs adaptivity to avoid a large dimensional coarse spaces. Because the coupling over the uncertainty space, the adaptivity
simply requires grouping relevant realizations.

We start with an arbitrary spatial coarse grid and perform an adaptivity in the uncertainty space by selecting
realizations. The latter requires some carefully selected distances as simple $L^2$ distance in the coefficients is an appropriate measure for multiscale problems. We then introduce an appropriate distance and present an approximation for the distance, which is used for clustering. Once the clustering is performed, we construct
multiscale basis functions for each cluster and perform multiscale simulations. We discuss some further extensions.
One of them includes the online procedures, which can be used to accelerate the convergence. Numerical results are presented to demonstrate the accuracy and efficiency of the proposed method for several multiscale stochastic problems without scale separation. Most importantly, we find that the cluster-based method allows us to explore the heterogeneities of the solution space and to solve nonlinear random elliptic PDEs.

\section*{Acknowledgements}

The research of E. Chung is supported by Hong Kong RGC General Research Fund (Project 14317516).
Y. Efendiev would like to thank the partial support from NSF 1620318, the U.S. Department of Energy Office of Science,
Office of Advanced Scientific Computing Research, Applied Mathematics program under Award
Number DE-FG02-13ER26165 and National Priorities Research Program grant NPRP grant 7-1482-1278 from the Qatar National Research Fund. The research of Z. Zhang is supported by the Hong Kong RGC grants (27300616, 17300817), National Natural Science Foundation of China (Project 11601457), and Seed Funding Programme for Basic Research (HKU).

\bibliographystyle{plain}
\bibliography{ZWpaper,references}

\begin{thebibliography}{10}

\bibitem{babuska:04}
I.~Babuska, R.~Tempone, and G.~Zouraris.
\newblock Galerkin finite element approximations of stochastic elliptic partial
  differential equations.
\newblock {\em SIAM J. Numer. Anal.}, 42:800--825, 2004.

\bibitem{Caflisch:98}
R.~E. Caflisch.
\newblock Monte {C}arlo and quasi-{M}onte {C}arlo methods.
\newblock {\em Acta Numerica}, 7:1--49, 1998.

\bibitem{randomized2014}
VM~Calo, Y~Efendiev, J~Galvis, and G~Li.
\newblock Randomized oversampling for generalized multiscale finite element
  methods.
\newblock {\em http://arxiv.org/pdf/1409.7114.pdf}.
\newblock to appear in SIAM MMS.

\bibitem{ChengHouYanZhang:13}
M.~Cheng, T.~Y. Hou, M.~Yan, and Z.~Zhang.
\newblock A data-driven stochastic method for elliptic {PDE}s with random
  coefficients.
\newblock {\em SIAM J. UQ}, 1:452--493, 2013.

\bibitem{ChengHouZhang1:13}
M.~Cheng, T.~Y. Hou, and Z.~Zhang.
\newblock A dynamically bi-orthogonal method for stochastic partial
  differential equations {I}: derivation and algorithms.
\newblock {\em J. Comput. Phys.}, 242:843--868, 2013.

\bibitem{ChengHouZhang2:13}
M.~Cheng, T.~Y. Hou, and Z.~Zhang.
\newblock A dynamically bi-orthogonal method for stochastic partial
  differential equations {II}: adaptivity and generalizations.
\newblock {\em J. Comput. Phys.}, 242:753--776, 2013.

\bibitem{chung2016adaptive}
Eric Chung, Yalchin Efendiev, and Thomas~Y Hou.
\newblock Adaptive multiscale model reduction with generalized multiscale
  finite element methods.
\newblock {\em Journal of Computational Physics}, 320:69--95, 2016.

\bibitem{chung2015residual}
Eric~T Chung, Yalchin Efendiev, and Wing~Tat Leung.
\newblock Residual-driven online generalized multiscale finite element methods.
\newblock {\em Journal of Computational Physics}, 302:176--190, 2015.

\bibitem{chung2017constraint}
Eric~T Chung, Yalchin Efendiev, and Wing~Tat Leung.
\newblock Constraint energy minimizing generalized multiscale finite element
  method.
\newblock {\em arXiv preprint arXiv:1704.03193}, 2017.

\bibitem{chung2014adaptive}
ET~Chung, Y~Efendiev, and G~Li.
\newblock An adaptive {GMsFEM} for high-contrast flow problems.
\newblock {\em Journal of Computational Physics}, 273:54--76, 2014.

\bibitem{ge09_2}
Y.~Efendiev and J.~Galvis.
\newblock A domain decomposition preconditioner for multiscale high-contrast
  problems.
\newblock In Y.~Huang, R.~Kornhuber, O.~Widlund, and J.~Xu, editors, {\em
  Domain Decomposition Methods in Science and Engineering XIX}, volume~78 of
  {\em Lect. Notes in Comput. Science and Eng.}, pages 189--196.
  Springer-Verlag, 2011.

\bibitem{egh12}
Y.~Efendiev, J.~Galvis, and T.~Hou.
\newblock Generalized multiscale finite element methods.
\newblock {\em Journal of Computational Physics}, 251:116--135, 2013.

\bibitem{Efendiev:13}
Y.~Efendiev, J.~Galvis, and T.~Y. Hou.
\newblock Generalized multiscale finite element methods.
\newblock {\em J. Comput. Phys.}, 251:116--135, 2013.

\bibitem{ehw99}
Y.~Efendiev, T.~Hou, and X.H. Wu.
\newblock Convergence of a nonconforming multiscale finite element method.
\newblock {\em SIAM J. Numer. Anal.}, 37:888--910, 2000.

\bibitem{HouEfendievWu:2000}
Y.~R. Efendiev, T.~Y. Hou, and X.~Wu.
\newblock Convergence of a nonconforming multiscale finite element method.
\newblock {\em SIAM Journal on Numerical Analysis}, 37(3):888--910, 2000.

\bibitem{Ghanem:91}
R.~G. Ghanem and P.~D. Spanos.
\newblock {\em Stochastic {F}inite {E}lements: a {S}pectral {A}pproach.}
\newblock Springer-Verlag, New York, 1991.

\bibitem{Giles:08}
M.~B. Giles.
\newblock Multilevel {M}onte {C}arlo path simulation.
\newblock {\em Operations Research}, 56:607--617, 2008.

\bibitem{glasserman:03}
P.~Glasserman.
\newblock {\em Monte {C}arlo methods in financial engineering}, volume~53.
\newblock Springer Science, 2003.

\bibitem{Tropp:11}
N.~Halko, P.~Martinsson, and J.~Tropp.
\newblock Finding structure with randomness: probabilistic algorithms for
  constructing approximate matrix decompositions.
\newblock {\em SIAM review}, 53:217--288, 2011.

\bibitem{Hartigan:1979}
J.A. Hartigan and M.~A. Wong.
\newblock Algorithm {A}{S} 136: {A} k-means clustering algorithm.
\newblock {\em Journal of the Royal Statistical Society. Series C (Applied
  Statistics)}, 28(1):100--108, 1979.

\bibitem{hw97}
T.~Hou and X.H. Wu.
\newblock A multiscale finite element method for elliptic problems in composite
  materials and porous media.
\newblock {\em J. Comput. Phys.}, 134:169--189, 1997.

\bibitem{WuanHou:06}
T.~Y. Hou, W.~Luo, B.~Rozovskii, and H.~M. Zhou.
\newblock Wiener chaos expansions and numerical solutions of randomly forced
  equations of fluid mechanics.
\newblock {\em J. Comput. Phys.}, 216:687--706, 2006.

\bibitem{Karhunen:47}
K.~Karhunen.
\newblock Uber lineare methoden in der {W}ahrscheinlichkeitsrechnung.
\newblock {\em Ann. Acad. Sci. Fennicae. Ser. A. I. Math.-Phys.}, 37:1--79,
  1947.

\bibitem{Rokhlin:07}
E.~Liberty, F.~Woolfe, P.G. Martinsson, V.~Rokhlin, and M.~Tygert.
\newblock Randomized algorithms for the low-rank approximation of matrices.
\newblock {\em Proceedings of the National Academy of Sciences},
  104(51):20167--20172, 2007.

\bibitem{Loeve:78}
M.~Lo\`{e}ve.
\newblock {\em Probability theory. {V}ol. II, 4th ed. {GTM}. 46.}
\newblock Springer-Verlag, ISBN 0-387-90262-7, 1978.

\bibitem{Knio:01}
O.~Ma{\i}tre, O.~M. Knio, H.~Najm, and R.~G. Ghanem.
\newblock A stochastic projection method for fluid flow: {I}. basic
  formulation.
\newblock {\em J. Comput. Phys.}, 173:481--511, 2001.

\bibitem{manning2008introduction}
Christopher~D Manning, Prabhakar Raghavan, Hinrich Sch{\"u}tze, et~al.
\newblock {\em Introduction to information retrieval}, volume~1.
\newblock Cambridge university press Cambridge, 2008.

\bibitem{matthies:05}
H.~G. Matthies and A.~Keese.
\newblock Galerkin methods for linear and nonlinear elliptic stochastic partial
  differential equations.
\newblock {\em Comput. Method Appl. Mech. Eng.}, 194:1295--1331, 2005.

\bibitem{Najm:09}
H.~N. Najm.
\newblock Uncertainty quantification and polynomial chaos techniques in
  computational fluid dynamics.
\newblock {\em Annual Review of Fluid Mechanics}, 41:35--52, 2009.

\bibitem{Niederreiter:01}
H.~Niederreiter.
\newblock {\em Quasi-{M}onte {C}arlo {M}ethods}.
\newblock Wiley Online Library., 2010.

\bibitem{nouy2007generalized}
Anthony Nouy.
\newblock A generalized spectral decomposition technique to solve a class of
  linear stochastic partial differential equations.
\newblock {\em Computer Methods in Applied Mechanics and Engineering},
  196(45):4521--4537, 2007.

\bibitem{nouy2010proper}
Anthony Nouy.
\newblock Proper generalized decompositions and separated representations for
  the numerical solution of high dimensional stochastic problems.
\newblock {\em Archives of Computational Methods in Engineering},
  17(4):403--434, 2010.

\bibitem{sapsis:09}
T.~Sapsis and P.~Lermusiaux.
\newblock Dynamically orthogonal field equations for continuous stochastic
  dynamical systems.
\newblock {\em Physica D: Nonlinear Phenomena}, 238:2347--2360, 2009.

\bibitem{Zabaras:13}
J.~Wan and N.~Zabaras.
\newblock A probabilistic graphical model approach to stochastic multiscale
  partial differential equations.
\newblock {\em Journal of Computational Physics}, 250:477--510, 2013.

\bibitem{Wan:06}
X.~L. Wan and G.~Karniadakis.
\newblock Multi-element generalized polynomial chaos for arbitrary probability
  measures.
\newblock {\em SIAM J. Sci. Comp.}, 28:901--928, 2006.

\bibitem{Xiu:09}
D.~Xiu.
\newblock Fast numerical methods for stochastic computations: a review.
\newblock {\em Commun. Comput. Phys.}, 5:242--272, 2009.

\bibitem{Xiu:03}
D.~Xiu and G.~Karniadakis.
\newblock Modeling uncertainty in flow simulations via generalized polynomial
  chaos.
\newblock {\em J. Comput. Phys.}, 187:137--167, 2003.

\bibitem{ZhangCiHouMMS:15}
Z.~Zhang, M.~Ci, and T.~Y. Hou.
\newblock A multiscale data-driven stochastic method for elliptic {PDE}s with
  random coefficients.
\newblock {\em SIAM Multiscale Model. Simul.}, 13:173--204, 2015.

\end{thebibliography}


\end{document}